\documentclass[twocolumns]{article}
\usepackage{epsfig}
\usepackage{graphicx}
\usepackage{color}
\usepackage[acmtitlespace,nocopyright]{acmneww}
\usepackage{times}
\usepackage{amsmath}
\usepackage{amssymb}
\usepackage{xspace}
\usepackage{txfonts}
\usepackage{enumitem}

\newcommand {\mm}[1] {\ifmmode{#1}\else{\mbox{\(#1\)}}\fi}

\newcommand{\proof}{\noindent{\sc Proof.~}}

\newcommand{\eop}{\hfill\usebox{\smallProofsym}\bigskip}  %
\newsavebox{\smallProofsym}                            % smallproofsymbol
\savebox{\smallProofsym}                               %
{%                                                     %
}                                                      %
\makeatletter
\long\def\@makecaption#1#2{%
  \vskip\abovecaptionskip
  \sbox\@tempboxa{\small #1: #2}%
  \ifdim \wd\@tempboxa >\hsize
    \small #1: #2\par
  \else
    \global \@minipagefalse
    \hb@xt@\hsize{\hfil\box\@tempboxa\hfil}%
  \fi
  \vskip\belowcaptionskip}
\makeatother

\newcommand{\Bspace}        {\mm{{\mathbb B}}}
\newcommand{\Rspace}        {\mm{{\mathbb R}}}

\newcommand{\Zspace}        {\mm{{\mathbb Z}}}

\newcommand{\Ecal}          {\mm{\cal E}}
\newcommand{\Fcal}          {\mm{\cal F}}
\newcommand{\Gcal}          {\mm{\cal G}}
\newcommand{\Scal}          {\mm{\cal S}}

\newcommand{\Delaunay}[1]   {\mm{\rm Del\,}{#1}}

\newcommand{\Area}[1]       {\mm{\rm Area}{({#1})}}
\newcommand{\Volume}[1]     {\mm{\rm Vol}{({#1})}}
\newcommand{\Simplex}       {\mm{\sigma}}
\newcommand{\Timplex}       {\mm{\tau}}
\newcommand{\Circumcenter}[1]  {\mm{\rm Circumcenter}{({#1})}}
\newcommand{\Circumradius}[2]  {\mm{\rm Circumradius}^{#1}{({#2})}}
\newcommand{\Inradius}[1]      {\mm{\rm Inradius}{({#1})}}
\newcommand{\Centroid}[1]      {\mm{\rm Centroid}{({#1})}}

\newcommand{\conv}[1]       {\mm{\rm conv\,}{#1}}
\newcommand{\Edist}[2]      {\mm{\|{#1}-{#2}\|}}
\newcommand{\norm}[1]       {\mm{\|{#1}\|}}

\newtheorem{result}{}

\newcommand{\Point}[1]       {}
\newcommand{\ignore}[1]{}

\title{Functionals on Triangulations of Delaunay Sets
       \thanks{This research is partially supported by
               the Russian Government under the Mega Project 11.G34.31.0053,
               RFBR grant 11-01-00735, DMS 1101688, and
               the European Science Foundation (ESF) under
                   the Research Network Programme.}
       }

\author{Nikolay P.\ Dolbilin\thanks{Steklov Mathematics Institute,
       Moscow, Russian Federation.},
        Herbert Edelsbrunner\thanks{IST Austria (Institute of Science and
          Technology Austria), Kloster\-neu\-burg, Austria,
          Departments of Computer Science and of Mathematics,
          Duke University, Durham, North Carolina,
          and Geomagic, Research Triangle Park, North Carolina.},
        Alexey Glazyrin\thanks{Mathematics Department,
          University of Texas at Brownsville, Texas, USA.}
        and Oleg R.\ Musin${}^{\S}$
}

\begin{document}
\maketitle

\begin{abstract}
  We study densities of functionals over uniformly bounded triangulations
  of a Delaunay set of vertices, and prove that the minimum
  is attained for the Delaunay triangulation if this is the case
  for finite sets.
\end{abstract}

\vspace{0.1in}

%%%%%%%%%%%%%%%%%%%%%%%%%%%%%%%%%%%%%%%%%%%%%%%%%%%%%%%%%%%%%%%%%%%%%%%%%%
%%%%%%%%%%%%%%%%%%%%%%%%%%%%%%%%%%%%%%%%%%%%%%%%%%%%%%%%%%%%%%%%%%%%%%%%%%
\section{Introduction}
\label{sec1}
%%%%%%%%%%%%%%%%%%%%%%%%%%%%%%%%%%%%%%%%%%%%%%%%%%%%%%%%%%%%%%%%%%%%%%%%%%
%%%%%%%%%%%%%%%%%%%%%%%%%%%%%%%%%%%%%%%%%%%%%%%%%%%%%%%%%%%%%%%%%%%%%%%%%%

A \emph{Delaunay set} $X \subseteq \Rspace^d$ has positive numbers
$r < R$ such that every open ball of radius $r$ contains at most one point,
and every closed ball of radius $R$ contains at least one point of $X$.
Such sets were introduced as $(r,R)$-systems by Boris N.\ Delaunay in 1924.
By a \emph{triangulation} of $X$, we mean a simplicial complex, $T$,
whose vertex set is $X$ and whose underlying space is $\Rspace^d$.
This triangulation is \emph{uniformly bounded} if there is
a real number $q = q(T)$ such that the circumsphere of every $d$-simplex
in $T$ has radius smaller than or equal to $q$.
A particular triangulation is the \emph{Delaunay triangulation},
denoted as $\Delaunay{X}$, 
whose $d$-simplices satisfy the additional condition that all other
vertices lie outside their circumspheres.
It exists if $X$ is generic, as will be explained shortly.
The Delaunay triangulation of a Delaunay set is necessarily uniformly bounded.
We also consider Delaunay triangulations of \emph{finite} sets of points,
for which the underlying space of the simplicial complex is the convex
hull of the points.

Writing $\Scal_d$ for the set of all $d$-simplices in $\Rspace^d$,
we consider functionals $F: \Scal_d \to \Rspace$ for which there
are constants $e = e(r,q,d)$ and $E = E(r,q,d)$ such that
$e \leq F(\Simplex) \leq E$ for every $d$-simplex $\Simplex$ whose
edges are longer than or equal to $2r$ and whose circumsphere
has a radius smaller than or equal to $q$.
Writing $\Ecal$ for this class of functionals,
we define subclasses $\Gcal \subseteq \Fcal \subseteq \Ecal$
by requiring additional conditions.
Briefly, $F$ belongs to $\Fcal$ if the sum of values over the $d$-simplices
of the Delaunay triangulation of $d+2$ points in $\Rspace^d$ is smaller
than or equal to the sum over the $d$-simplices in the other triangulation,
and $F$ belongs to $\Gcal$ is it satisfies a similar condition
for all finite sets of points.
For a triangulation, we define the \emph{density} of the functional
by taking sums and lower limits over a growing sequence of balls:
\begin{eqnarray}
  f(T)  &=&  \liminf_{\alpha \to \infty} \frac{1}{\Volume{\Bspace_\alpha}}
             \sum_{\Bspace_\alpha \supseteq \Simplex \in T} F(\Simplex) ,
\end{eqnarray}
where $\Bspace_\alpha$ is the closed ball with radius $\alpha$ and
center at the origin of $\Rspace^d$, and $\Volume{\Bspace_\alpha}$
is its volume.
With these definitions, we can give our main result:
\begin{itemize}
  \item  in $\Rspace^d$, $F \in \Gcal$ implies that the Delaunay triangulation
    minimizes the density of $F$ among all uniformly bounded triangulations
    of a Delaunay set,
    and in $\Rspace^2$, $F \in \Fcal$ suffices to reach the same conclusion.
\end{itemize}
There are many concrete functionals studied in the literature
to which our result applies.
Here, we just mention two:
\begin{itemize}
  \item  the functional that maps every triangle in $\Rspace^2$
    to the radius of its circumcircle; see \cite{Mus97},
  \item  the functional that maps every $d$-simplex to the sum
    of squares of its edge lengths times the volume; see \cite{Raj94}.
\end{itemize}
The remainder of this paper presents the detailed results in
two sections.

%\clearpage
%%%%%%%%%%%%%%%%%%%%%%%%%%%%%%%%%%%%%%%%%%%%%%%%%%%%%%%%%%%%%%%%%%%%%%%%%%
%%%%%%%%%%%%%%%%%%%%%%%%%%%%%%%%%%%%%%%%%%%%%%%%%%%%%%%%%%%%%%%%%%%%%%%%%%
\section{Background}
\label{sec2}
%%%%%%%%%%%%%%%%%%%%%%%%%%%%%%%%%%%%%%%%%%%%%%%%%%%%%%%%%%%%%%%%%%%%%%%%%%
%%%%%%%%%%%%%%%%%%%%%%%%%%%%%%%%%%%%%%%%%%%%%%%%%%%%%%%%%%%%%%%%%%%%%%%%%%

In this section, we introduce the background on Delaunay sets,
their uniformly bounded triangulations,
and functionals on such triangulations.

%%%%%%%%%%%%%%%%%%%%%%%%%%%%%%%%%%%%%%%%%%%%%%%%%%%%%%%%%%%%%%%%%%%%%%%%%%
\subsection{Delaunay Sets}
\label{sec21}
%%%%%%%%%%%%%%%%%%%%%%%%%%%%%%%%%%%%%%%%%%%%%%%%%%%%%%%%%%%%%%%%%%%%%%%%%%

We recall from Section \ref{sec1} that $X \subseteq \Rspace^d$ is a
\emph{Delaunay set} if there are positive constants $r < R$ such that
(I) every open ball of radius $r$ contains at most one point of $X$, and
(II) every closed ball of radius $R$ contains at least one point of $X$.
Hence, $X$ has no tight cluster and leaves no large hole.

\paragraph{Counting points.}
Condition (I) implies that every bounded subset of $\Rspace^d$ contains
only finitely many points of $X$.
Indeed, the subset can be covered by finitely many open balls of
radius $r$, and each such ball contains at most one point.
Condition (II) implies that every cone with non-zero volume
contains infinitely many points of $X$.
Indeed, the cone contains an infinite string of disjoint closed balls
of radius $R$, and each such ball contains at least one point of $X$.
We quantify the first observation by giving concrete estimates.
Let $\Bspace_\alpha (z)$ be the closed ball with radius $\alpha$
and center $z$, and call the difference between two concentric
balls an \emph{annulus}.
\begin{result}[Point Count Lemma]
  Let $X$ be a Delaunay sets with parameters $r < R$ in $\Rspace^d$.
  \begin{enumerate}
    \item[(i)] There are constants $p = p(R,d)$ and $P = P(r,d)$ such that the
      number of points of $X$ in $\Bspace_\alpha (z)$ is between
      $p \alpha^d$ and $P \alpha^d$.  \Point{Lemma 3.2.}
    \item[(ii)] There is a constant $P' = P' (r,d)$
      such that the number of points of $X$ in
      $\Bspace_{\alpha+1} (z) - \Bspace_\alpha (z)$ is at most
      $P' \alpha^{d-1}$.  \Point{Lemma 3.3.}
  \end{enumerate}
\end{result}
\proof
 To prove (i), we note that $\Bspace_\alpha (z)$ can be covered by some
 constant times $( {\alpha}/{r} )^d$ balls of radius $r$,
 and that we can pack some other constant times 
 $( {\alpha}/{R} )^d$ balls of radius $R$ in it.
 The lower and upper bounds follow.

 To prove (ii), we cover the annulus with a constant times
 $\alpha^{d-1}/r^d$ balls of radius $r$.
 The upper bound follows.
\eop

It should be clear that the bound in (ii) also holds for annuli
of constant width, but not for annuli whose width is a positive fraction
of the radius.

\paragraph{Delaunay triangulations.}
Following the original idea of Boris N.\ Delaunay, we consider $d$-simplices
with vertices from $X$ such that the open ball bounded by the
$(d-1)$-dimensional circumsphere contains no points of $X$.
We call such $d$-simplices \emph{empty}.
Here, it is convenient to assume that $X$ is \emph{generic} in the sense
that no $d+2$ points in $X$ lie on a common $(d-1)$-sphere.
Under this assumption, the empty $d$-simplices fit together without gap
and overlap.
Now consider two not necessarily empty but non-overlapping $d$-simplices
that share a $(d-1)$-simplex, which is a face of both.
Assuming the two $d$-simplices belong to a triangulation,
we call this face \emph{locally Delaunay} if the $(d+1)$-st vertex of
the second $d$-simplex lies outside the circumsphere of the first $d$-simplex.
Note that the condition is symmetric because the two circumspheres
intersect in the $(d-2)$-sphere that passes through the vertices of the face,
and either both $(d+1)$-st vertices lie outside or both lie inside
the respective other circumsphere.
Delaunay considered both conditions and proved that they are
equivalent \cite{Del34}.
\begin{result}[Delaunay Triangulation Theorem]
  Let $X$ be a generic Delaunay set in $\Rspace^d$.
  \begin{enumerate}
    \item[(i)] The collection of empty $d$-simplices together with their faces
      form a triangulation of $X$,
      commonly known as the \emph{Delaunay triangulation}, $\Delaunay{X}$.
      \Point{Theorem 1.3.}
    \item[(ii)] If all $(d-1)$-simplices of a triangulation $T$ of $X$
      are locally Delaunay, then $T = \Delaunay{X}$.
      \Point{Lemma 1.4.}
  \end{enumerate}
\end{result}
The equivalence between the local and the global conditions expressed
in (ii) also holds for finite sets $X$.
In the plane, it means that a triangulation of a generic set $X$ is Delaunay
iff for each edge the sum of opposite angles in the two incident triangles
is less than $\pi$.

%%%%%%%%%%%%%%%%%%%%%%%%%%%%%%%%%%%%%%%%%%%%%%%%%%%%%%%%%%%%%%%%%%%%%%%%%%
\subsection{Uniformly Bounded Triangulations}
\label{sec22}
%%%%%%%%%%%%%%%%%%%%%%%%%%%%%%%%%%%%%%%%%%%%%%%%%%%%%%%%%%%%%%%%%%%%%%%%%%

Let $X$ be a generic Delaunay set in $\Rspace^d$, and let $T$ be
a triangulation of $X$.
We recall that this means that $T$ is a
simplicial complex with vertex set $X$ whose underlying space is $\Rspace^d$.
Recall also that $T$ is \emph{uniformly bounded} if there is a real number
$q = q(T)$ such that the radius of the circumsphere of every $d$-simplex
in $T$ is smaller than or equal to $q$.
It follows that no edge of $T$ is longer than $2 q$.
Note that the Delaunay triangulation of $X$ is uniformly bounded
with $q = R$.

\paragraph{Not every triangulation is uniformly bounded.}
We begin by showing that every Delaunay set has triangulations that
are not uniformly bounded.
Given $X$, we construct such a triangulation in three steps.
\begin{enumerate}
  \item[1.]  For every point $x \in X$ and every $L > 0$, we can find
    many points $y \in X$ such that the edge $xy$ is longer than $L$
    and does not pass through any other points of $X$.
    Indeed, there is such an edge near every direction out of $x$.
    \Point{Lemma 2.1.}
    To see this, we consider the set of points in $X$ that lie within
    the closed ball of radius $L$ around $x$.
    There are only finitely many such points, which implies that
    within each cone with non-zero volume and apex $x$, we can find a 
    subcone, again with non-zero volume and apex $x$,
    that does not contain any of the points of $X$ inside the ball.
    However, as argued in Section \ref{sec21}, the cone contains
    infinitely many points of $X$, so they must all be at
    distance larger than $L$ from $x$.
    Among these, let $y$ be the point closest to $x$.
\end{enumerate}
Using the knowledge about long edges,
we construct a triangulation inductively, one phase at a time.
After the $k$-th phase, we will have a triangulation $T_k$ of a finite
subset of $X$ that includes all points at distance $k$ or less from
the origin.
In addition, we will make sure that $T_k$ contains at least one edge
longer than $k$, and that every point of $X$ that belongs to the
underlying space of $T_k$ is a vertex of $T_k$.
We start with $T_0$ consisting of a single edge connecting the
two points of $X$ that are closest to the origin of $\Rspace^d$.
\begin{enumerate}
  \item[2.]  In the $(k+1)$-st phase, we let $x$ be a vertex in
    the boundary of $T_k$.
    Let $y$ be another point of $X$ such that the edge $xy$ is longer
    than $k+1$ and does not intersect the simplices in $T_k$ other than at $x$.
    Since $T_k$ is a triangulation, its underlying space is convex,
    and its boundary is triangulated.
    Let $\Simplex$ be an $i$-simplex in the boundary of $T_k$ that is
    visible from $y$.
    We extend $T_k$ by adding the $(i+1)$-simplex formed by $y$
    and the vertices of the $i$-simplex.
    Doing this for $y$ and all visible simplices in the boundary of $T_k$,
    we obtain a simplicial complex $T_k'$ by \emph{starring} from $y$.
    Similarly, we add a point $z$ at distance $k+1$ or less from the
    origin that lies outside the underlying space by starring to the
    triangulation.
    Repeating this operation for all such points $z$,
    we eventually get a simplicial complex $T_k''$.
    \Point{Theorem 2.2, first half.}
  \item[3.]
    While $T_k''$ is a valid simplicial complex in $\Rspace^d$,
    some of the new simplices may contain points of $X$ in their interiors.
    By construction, $xy$ is not among these simplices.
    Let $w \in X$ be such a point, and $\Simplex \in T_k''$ the simplex
    of lowest dimension, $j$, that contains $w$ in its interior.
    Note that $j \geq 1$.
    We fix the situation by decomposing $\Simplex$ into $j+1$
    $j$-simplices, each the convex hull of $w$ and $j$ vertices of $\Simplex$.
    Similarly, we decompose each simplex that contains $\Simplex$ as
    a face into $j+1$ simplices of the same dimension by starring from $w$.
    Repeating this procedure for all such points $w$, we eventually
    get a triangulation $T_{k+1}$ such that all points in $X$
    that belong to the underlying space of $T_{k+1}$
    are in fact vertices of $T_{k+1}$.
    \Point{Theorem 2.2, second half.}
\end{enumerate}
Observe that the edge $xy$ added to $T_k$ in Step 2 remains undivided
until the end.
This implies that $T_{k+1}$ indeed contains an edge longer than $k+1$,
as required.
It follows that the triangulation thus constructed by transfinite induction
is not uniformly bounded.

\paragraph{Measuring volume.}
The remainder of this section states and proves properties of
uniformly bounded triangulations.
We begin with the volume of their simplices.
\begin{result}[Volume Lemma]
  Let $X$ be a Delaunay set with parameters $r < R$ in $\Rspace^d$,
  and let $T$ be a uniformly bounded triangulation with parameter $q$ of $X$.
  \begin{enumerate}
    \item[(i)] In $\Rspace^2$, there is a positive constant $v = v(r,q)$
      such that $v \leq \Area{\Simplex}$ for every triangle $\Simplex$ in $T$.
    \item[(ii)] In $\Rspace^d$, there is a constant $V = V(q,d)$ such that
      $\Volume{\Simplex} \leq V$ for every $d$-simplex $\Simplex$ in $T$.
      \Point{Lemma 2.3.}
  \end{enumerate}
\end{result}
\proof
 We prove (i) by expressing the area of a triangle
 in terms of the three edge lengths and the radius of the circumcircle:
 $\Area{\Simplex} = \frac{abc}{4 \varrho}$.
 The edges cannot be shorter than $2r$,
 and the radius cannot be larger than $q$, which implies
 $v(r,q) = 2r^3/q \leq \Area{\Simplex}$.

 To prove (ii), we note that every $d$-simplex is contained in
 the ball bounded by its circumsphere.
 Since the radius is at most $q$, this ball is smaller than
 $V(q,d) = (2q)^d > \Volume{\Simplex}$.
\eop

If we remove the requirement of uniform boundedness, then the proof
of the Volume Lemma breaks down.
It is not clear whether the upper bound fails.
In this context, we mention a related question asked by
L.\ Danzer and independently by M.\ Boshernitzen:
``{\it is it true that for every planar Delaunay set there exists
a triangle with arbitrarily large area that contains no points in
its interior?}''
This question is still open.

\begin{figure}[hbt]
 \centering
 \resizebox{!}{1.5in}{\input{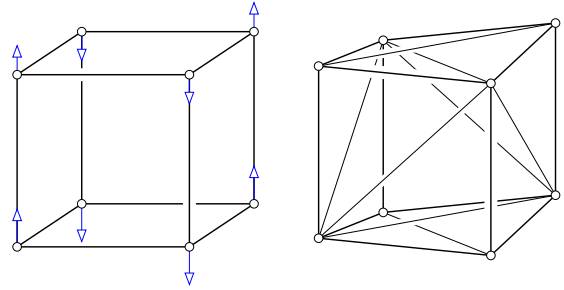}}
 \caption{A distorted cube decomposed into seven tetrahedra.
   The arrows indicating the distortion are exaggerated for
   better visibility.}
\label{fig:cube}
\end{figure}

Next, we describe a Delaunay set in $\Rspace^3$ that has tetrahedra
of arbitrarily small volume in the Delaunay triangulation.
It shows that the limitation of the lower bound in (i) to two dimensions
is necessary.
Consider the standard cubic lattice, $\Zspace^3$.
Let $\delta_\ell = \frac{1}{2 + |\ell|}$, for every $\ell \in \Zspace$,
and move every point $(i, j, k) \in \Zspace^3$ to
$(i, j, k+(-1)^{i+j} \delta_k)$, denoting the new point set by $X$.
To study the volume of the tetrahedra in the Delaunay triangulation,
we consider a single integer cube,
for which we get a tetrahedron of volume about $\frac{1}{3}$ in the middle,
four tetrahedra of volume about $\frac{1}{6}$ across each face,
and two flat tetrahedra at the top and the bottom;
see Figure \ref{fig:cube}.
These volume estimates assume arbitrarily small values of $\delta_\ell$.
If the third coordinates of the original vertices are $k$ and $k+1$,
then the volume of the top tetrahedron is $\frac{2}{3} \delta_{k+1}$,
and that of the bottom tetrahedron is $\frac{2}{3} \delta_k$.
Since among the $\delta_\ell$ there are arbitrarily small numbers,
there are tetrahedra in $\Delaunay{X}$ whose volume is arbitrarily
close to $0$.
It is not difficult to extend this example to four and higher dimensions.
\Point{Theorem 2.4.}

\paragraph{Counting simplices.}
Recall the Volume Lemma, which states that every triangle in a
uniformly bounded triangulation of a Delaunay set in $\Rspace^2$ 
has an area that exceeds a positive constant.
Since a disk of radius $\alpha$ has area $\alpha^2 \pi$,
this implies that the number of triangles contained in this disk
is at most some constant times $\alpha^2$.
A similar result holds in three and higher dimensions,
but the lack of a lower bound on the volume of a $d$-simplex
requires a different argument, which we present as the proof
of the following bounds.
\begin{result}[Simplex Count Lemma]
  Let $X$ be a Delaunay set with parameters $r < R$ in $\Rspace^d$,
  and let $T$ be a uniformly bounded triangulation with parameter $q$ of $X$.
  \begin{enumerate}
    \item[(i)] The number of $d$-simplices sharing a common vertex
      is bounded from above by a constant $S'' = S'' (r, q, d)$.
      \Point{Lemma 3.4.}
    \item[(ii)] There are positive constants $s = s(R, q, d)$
      and $S = S(r, q, d)$ such that the number of simplices contained
      in a ball of radius $\alpha > 4q$ is between
      $s \alpha^d$ and $S \alpha^d$.  \Point{Lemma 3.5.}
    \item[(iii)] There is a constant $S' = S'(r, q, d)$ such that
      the number of $d$-simplices contained in a ball of radius $\alpha+1$
      but not in the concentric ball of radius $\alpha$ is at most
      $S' \alpha^{d-1}$.
      \Point{Lemma 3.6.}
  \end{enumerate}
\end{result}
\proof
 To prove (i), we let $x$ be the shared vertex, and we note that all
 incident $d$-simplices are contained in the ball of radius $2q$
 centered at $x$.
 By the Point Count Lemma, the number of points in this ball
 is bounded from above by $P(r, d) \cdot (2q)^d$.
 We have at most one $d$-simplex for every combination of $d$ of
 these points, which gives $S'' (r, q, d)$.

 The upper bound in (ii) is now easy:
 by the Point Count Lemma, the number of points inside the ball of
 radius $\alpha$ is at most $P(r, d) \cdot \alpha^d$.
 Multiplying $P(r, d)$ with $S'' (r, q, d)$ gives $S (r, q, d)$.
 To get the lower bound, we restrict ourselves to the ball of
 radius $\alpha - 2q$.
 By the Point Count Lemma, the number of points in this smaller
 ball is at least $p(R, q) \cdot (\alpha - 2q)^d$.
 Every $d$-simplex incident to one of these points is contained
 in the ball of radius $\alpha$.
 Each point belongs to at least $d+1$ $d$-simplices, which implies
 that the lower bound on the number of points also applies to the
 $d$-simplices.
 Finally, $(\alpha - 2q)^d$ is at least $\alpha^d / 2^d$.

 To prove (iii), we use the upper bound of $P' (r, d) \cdot \alpha^{d-1}$
 on the number of points in the annulus.
 Each $d$-simplex we count is incident to at least one of these points.
 Multiplying $P' (r, d)$ with $S''(r, q, d)$ gives $S'(r, q, d)$.
\eop

%%%%%%%%%%%%%%%%%%%%%%%%%%%%%%%%%%%%%%%%%%%%%%%%%%%%%%%%%%%%%%%%%%%%%%%%%%
\subsection{Functionals}
\label{sec23}
%%%%%%%%%%%%%%%%%%%%%%%%%%%%%%%%%%%%%%%%%%%%%%%%%%%%%%%%%%%%%%%%%%%%%%%%%%

Recall that $\Scal_d$ denotes the set of simplices,
including degenerate ones.
We are interested in functionals that have constant
upper and lower bounds for the simplices that arise in
uniformly bounded triangulations of Delaunay sets.
For other degenerate simplices we also allow infinity as a value.
\begin{result}[Definition]
  Let $\Ecal$ be the class of functionals $F: \Scal_d \to \Rspace$
  for which there are constants $e = e(r,q,d)$ and $E = E(r,q,d)$
  such that $e \leq F(\Simplex) \leq E$ for all $d$-simplices
  $\Simplex$ with edges of length at least $2r$ and radius of the
  circumsphere at most $q$.
\end{result}
\ignore{
  The space of such simplices is compact, which implies that a functional
  that is continuous in the edge lengths belongs to $\Ecal$.}
In this section, we extend the functionals from simplices to
triangulations, and we introduce subclasses that favor
Delaunay triangulations for finite sets of points.

\paragraph{Densities.}
As already mentioned in Section \ref{sec1}, we define the density
of a functional on a triangulation by taking the lower limit
over a growing ball,
of the sum of values over all $d$-simplices in the ball divided by
the volume of the ball:
\begin{eqnarray}
  f(T)  &=&  \liminf_{\alpha \to \infty} \frac{1}{\Volume{\Bspace_\alpha}}
             \sum_{\Bspace_\alpha \supseteq \Simplex \in T} F(\Simplex) .
  \label{eqn:functional}
\end{eqnarray}
There are other possibilities, such as taking the upper limit,
or taking the average over the simplices.
Our results extend to both modifications of the definition.
One of Delaunay's motivations for defining $(r, R)$-systems
was to generalize lattices in $\Rspace^d$ to a larger class of sets.
For the Delaunay triangulation of any lattice $\Lambda \subseteq \Rspace^d$,
the limit of the expression in \eqref{eqn:functional},
in which we substitute $\Delaunay{\Lambda}$ for $T$, is well defined.
Unfortunately, this is not generally the case for Delaunay triangulations
of Delaunay sets, which is the reason for taking the lower limit.
Since this might not be entirely obvious, we will prove shortly
that for a broad class of functionals in $\Ecal$, the limit does not
generally exist.
Before that, we prove some positive results,
namely that the density of every functional
is bounded and independent of the choice of origin.
Specifically, we define
\begin{eqnarray}
  f_z(T)  &=&  \liminf_{\alpha \to \infty} \frac{1}{\Volume{\Bspace_\alpha (z)}}
             \sum_{\Bspace_\alpha (z) \supseteq \Simplex \in T} F(\Simplex) 
  \label{eqn:functional-z}
\end{eqnarray}
for every point $z \in \Rspace^d$,
and we prove that all choices of $z$ give the same result.
\begin{result}[Properties]
  Let $F$ be a functional in $\Ecal$.
  \begin{enumerate}
    \item[(i)] There is a constant $C = C(r, q, d)$ such that $f(T) \leq C$
      for every uniformly bounded triangulation of a Delaunay set
      in $\Rspace^d$.
      \Point{Corollary 3.7.}
    \item[(ii)] $f(T) = f_z (T)$ for every $z \in \Rspace^d$. 
      \Point{Lemma 3.8.}
  \end{enumerate}
\end{result}
\proof
 To prove (i), we recall the Simplex Count Lemma,
 which implies that the number of $d$-simplices
 contained in $\Bspace_\alpha$ is bounded from above by $S (r, q, d)$.
 Multiplying with $E(r,q,d)$ gives $C(r, q, d)$.

 To prove (ii), we let $f(T,\alpha)$ and $f_z(T,\alpha)$ be the expressions
 in \eqref{eqn:functional} and \eqref{eqn:functional-z}
 without taking the lower limit,
 so that $f(T) = \liminf_{\alpha \to \infty} f(T,\alpha)$,
 and similarly for $f_z(T)$ and $f_z(T,\alpha)$.
 It suffices to prove
 \begin{eqnarray}
   \lim_{\alpha \to \infty} [f(T,\alpha)-f_z(T,\alpha)]  &=&  0 ,
   \label{eqn:limit1}
 \end{eqnarray}
 which we do in two steps, namely by proving
 \begin{eqnarray}
   \lim_{\alpha \to \infty} [f(T,\alpha+L)-f(T,\alpha)]  &=&  0 ,
     \label{eqn:limit2} \\
   \lim_{\alpha \to \infty} [f(T,\alpha+L)-f_z(T,\alpha)]  &=&  0 ,
     \label{eqn:limit3}
 \end{eqnarray}
 where $L = \norm{z}$ is the distance between $z$ and the origin.
 To prove \eqref{eqn:limit2}, we write $V_d$ for the volume of the
 unit ball in $\Rspace^d$,
 and we note that $\Volume{\Bspace_\alpha} = V_d \alpha^d$.
 Furthermore, we write $\Sigma_B$ and $\Sigma_A$ for the sums over the
 $d$-simplices contained in the smaller ball and the extra $d$-simplices
 contained in the larger ball:
 \begin{eqnarray}
   \Sigma_B        &=&  \sum_{\Bspace_\alpha \supseteq \Simplex \in T} F(\Simplex) ,
     \label{eqn:SB} \\
   \Sigma_A + \Sigma_B  &=&  \sum_{\Bspace_{\alpha + L} \supseteq \Simplex \in T} F(\Simplex) .
     \label{eqn:SA}
 \end{eqnarray}
 By the Simplex Count Lemma, $\Sigma_B$ is at most some positive constant
 times $\alpha^d$, while $\Sigma_A$ is at most some constant times $\alpha^{d-1}$.
 Hence,
 \begin{eqnarray}
   \Delta  &=&  f(T,\alpha+L) - f(T,\alpha)                                 \\
           &=&  \frac{\Sigma_B+\Sigma_A}{V_d (\alpha+L)^d} - \frac{\Sigma_B}{V_d \alpha^d} \\
           &=&  \frac{\Sigma_A}{V_d (\alpha+L)^d}
             - \frac{(\alpha+L)^d \Sigma_B - \alpha^d \Sigma_B}
                         {V_d (\alpha+L)^d \alpha^d} .
   \label{eqn:Delta}
 \end{eqnarray}
 The first term in \eqref{eqn:Delta} goes to zero because
 $\Sigma_A$ grows slower than $\alpha^d$,
 and the second term goes to zero because
 $(\alpha+L)^d - \alpha^d$ grows slower than $\alpha^d$.
 The argument for \eqref{eqn:limit3} is similar.
 Indeed, all we need is to notice that the set of $d$-simplices
 contained in $\Bspace_{\alpha+L}$ but not contained in $\Bspace_\alpha (z)$
 is a subset of those contained in $\Bspace_{\alpha+L}$ but not
 contained in $\Bspace_{\alpha-L}$.
 By the Simplex Count Lemma, the number of simplices thus defined
 is bounded from above by a constant times $\alpha^{d-1}$,
 so that the argument goes through as before.
\eop
 
\paragraph{Non-existence of limits.}
We now show that taking the lower limit in the definition of density
is necessary because the limit does not generally exist.
Indeed, functionals for which the limit exists, even just for all
Delaunay triangulations of Delaunay sets, are the exception.
This is true in particular for the functionals that are invariant
under isometries, which include all examples we discuss in this paper.

\begin{figure}[hbt]
 \centering
 \resizebox{!}{1.7in}{\input{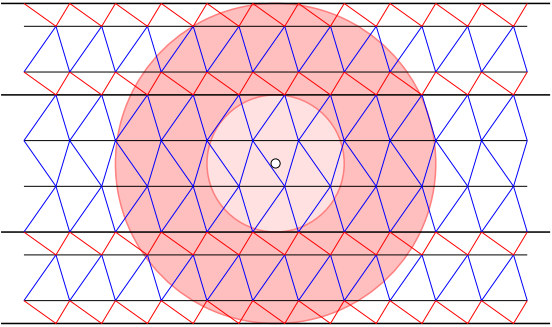}}
 \caption{A Delaunay triangulation made of a block of three wide strips
   in the center and two blocks alternating between narrow and wide strips
   glues on both sides.}
 \label{fig:strips}
\end{figure}

We begin by exhibiting a construction in $\Rspace^2$ that acts as a
stepping stone in our argument.
Let $\Simplex$ and $\Timplex$ be two \emph{compatible} triangles,
by which we mean that they share an edge, the two angles opposite
that edge add up to less than $\pi$, and the remaining four angles
are all acute.
The condition implies that at least one of the triangles is acute,
and we assume $\Simplex$ is.
Using a linear sequence of congruent copies of $\Simplex$,
we form a strip $T_\Simplex$, which we call \emph{wide},
and using copies of $\Timplex$,
we form a strip $T_\Timplex$, which we call \emph{narrow};
see Figure \ref{fig:strips}.
Gluing strips together so that they match up at boundary edges,
we get a Delaunay triangulation, provided no two narrow strips are glued
to each other.
Let $m_1, m_2, \ldots$ be an infinite sequence of odd integers.
We construct a Delaunay triangulation $D$ inductively,
starting with a block of $m_1$ wide strips.
On each side, we add a block of $m_2$ strips alternating between
narrow and wide,
then a block of $m_3$ wide strips,
then a block of $m_4$ strips again alternating between narrow and wide,
and so on.
For each $i \geq 1$, let $2 \alpha_i$ be the total width of
the first $2i-1$ blocks.
Making sure that the origin lies on the center line of the first block,
$B_{\alpha_i}$ is the largest disk centered at the origin that is still
contained in the union of the first $2i-1$ blocks.
We consider the sequence
\begin{eqnarray}
  f_i  &=&  \frac{1}{\alpha_i^2 \pi}
            \sum_{\Bspace_{\alpha_i} \supseteq \Simplex \in D} F(\Simplex) .
\end{eqnarray}
Assuming that the limit in \eqref{eqn:functional} exists,
the sequence of $f_i$ must converge for every sequence of $m_i$.
This is indeed the case if we measure area, because the number of
triangles that intersect $\Bspace_\alpha$ but are not contained in it
is bounded from above by a constant times $\alpha$.
Hence, $\lim_{\alpha \to \infty} \frac{1}{\alpha^2 \pi} \sum A(\Simplex) = 1$,
where we abbreviate $A(\Simplex) = \Area{\Simplex}$ and take the sum
over all triangles contained in $\Bspace_\alpha$, as usual.
This motivates us to consider the ratios of the terms in the two sequences.
Assuming $\Bspace_{\alpha_i}$ contains $k_i$ congruent copies of $\Simplex$
and $\ell_i$ congruent copies of $\Timplex$, this gives
\begin{eqnarray}
  g_i  &=&  \frac{\sum_{\Bspace_{\alpha_i} \supseteq \Simplex \in D} F(\Simplex)}
                 {\sum_{\Bspace_{\alpha_i} \supseteq \Simplex \in D} A(\Simplex)}
      ~~=~~  \frac{k_i F(\Simplex) + \ell_i F(\Timplex)}
                  {k_i A(\Simplex) + \ell_i A(\Timplex)} .
\end{eqnarray}
Define $Q_\Simplex = \frac{F(\Simplex)}{A(\Simplex)}$,
$Q_\Timplex = \frac{F(\Timplex)}{A(\Timplex)}$, and
$Q = \frac{F(\Simplex) + F(\Timplex)}{A(\Simplex) + A(\Timplex)}$.
Assuming $Q_\Simplex \neq Q_\Timplex$, we have
$Q_\Simplex \neq Q$ and define $\Delta = | Q_\Simplex - Q |$.
By definition, $g_1 = Q_\Simplex$.
We choose $m_2$ large enough so that
$|g_2 - Q| < \frac{\Delta}{3}$,
which is possible because $\ell_2 / k_2$ goes to $1$ as $m_2$ goes to infinity.
Then we choose $m_3$ large enough so that
$|g_3 - Q_\Simplex| < \frac{\Delta}{3}$,
and so on,
alternating between being close to $Q_{\Simplex}$ and $Q$.
We thus arrive at a contradiction because there is a gap of size
$\frac{\Delta}{3}$ between the terms with odd and even indices.
In other words, we need
$\frac{F(\Simplex)}{A(\Simplex)} = \frac{F(\Timplex)}{A(\Timplex)}$
for the limit to exist.

We finally show that the non-existence of the limit is not an artifact
of the particular Delaunay set we used in the construction of $D$.
Let $\Simplex'$ and $\Timplex'$ be arbitrary triangles
with longest edges of lengths $a$ and $c$ and opposite angles
$2 \varphi$ and $2 \psi$.
Setting $L > \max \{ \frac{a}{2 \cos \varphi}, \frac{c}{2 \cos \psi} \}$,
we construct triangles $\Simplex$ with edges of length $L, L, a$
and $\Timplex$ with edges of length $L, L, c$.
It is easy to verify that $\Simplex'$ and $\Simplex$ are compatible,
and so are $\Simplex$ and $\Timplex$, and $\Timplex$ and $\Timplex'$.
Repeating the construction with the strips and blocks three times,
we see that if $F$ is invariant under isometries of $\Rspace^2$
and the limit exists for the Delaunay triangulations of all Delaunay sets,
then $\frac{F(\Simplex')}{A(\Simplex')} = \frac{F(\Timplex')}{A(\Timplex')}$.
Conversely, among the functionals invariant under isometries,
only the ones proportional to the area have the limit defined
for the Delaunay triangulations of all Delaunay sets.
\Point{Theorem 3.1.}

\paragraph{Subclasses.}
We are interested in two subclasses of functionals,
$\Gcal \subseteq \Fcal \subseteq \Ecal$, which we now introduce.
To define $\Fcal$, let $Y$ be a generic set of $d+2$ points in $\Rspace^d$
such that no point lies inside the convex hull of the others.
The non-degenerate $d$-simplices spanned by the points cover the
convex hull twice; see Radon \cite{Rad21}.
Indeed, we can split them into two collections such that each forms
a triangulation of $Y$:
the Delaunay triangulation, $D = \Delaunay{Y}$,
and the other triangulation, $T$.
Changing one triangulation into the other is a \emph{flip},
a name motivated by the planar case in which it replaces
one diagonal of a convex quadrilateral with the other.
We give the flip a direction, leading from $T$ to $D$.
Let now $F$ be a functional, let $\Sigma_T$ be the sums of $F(\Simplex)$
over all $d$-simplices in $T$, and define $\Sigma_D$ similarly.
\begin{result}[Definition]
  The class $\Fcal$ consists of all functionals $F \in \Ecal$
  for which $\Sigma_D \leq \Sigma_T$.
\end{result}
In $\Rspace^2$, the extra property of functionals in $\Fcal$ suffices
to prove our main result.
In $\Rspace^d$, for $d \geq 3$, we need more structure.
The reason is the existence of triangulations that cannot be turned
into the Delaunay triangulation by a sequence of directed flips;
see \cite{Joe89} for finite examples in $\Rspace^3$.
Such examples do not exist in $\Rspace^2$; see \cite{Law77}.

Let now $Y$ be a finite set of points in $\Rspace^d$.
As before, we assume that $Y$ is generic.
Let $T'$ be a simplicial complex with vertex set $Y$,
but note that we do not require that $T'$ be a triangulation of $Y$.
For example, we could start with a triangulation of $Y$
and construct $T'$ as the subset of $d$-simplices that do not belong
to the Delaunay triangulation together with their faces.
Let $D'$ be the subset of simplices in $\Delaunay{Y}$
contained in the underlying space of $T'$.
Finally, let $\Sigma_{T'}$ be the sum of $F(\Simplex)$ over all
$d$-simplices in $T'$, and define $\Sigma_{D'}$ similarly.
\begin{result}[Definition]
  The class $\Gcal$ consists of all functionals $F \in \Ecal$
  for which $\Sigma_{D'} \leq \Sigma_{T'}$.
\end{result}
The condition for $F$ to belong to $\Gcal$ is at least as strong as
that for $F$ to belong to $\Fcal$, which implies $\Gcal \subseteq \Fcal$.

%\clearpage
%%%%%%%%%%%%%%%%%%%%%%%%%%%%%%%%%%%%%%%%%%%%%%%%%%%%%%%%%%%%%%%%%%%%%%%%%%
%%%%%%%%%%%%%%%%%%%%%%%%%%%%%%%%%%%%%%%%%%%%%%%%%%%%%%%%%%%%%%%%%%%%%%%%%%
\section{Results}
\label{sec3}
%%%%%%%%%%%%%%%%%%%%%%%%%%%%%%%%%%%%%%%%%%%%%%%%%%%%%%%%%%%%%%%%%%%%%%%%%%
%%%%%%%%%%%%%%%%%%%%%%%%%%%%%%%%%%%%%%%%%%%%%%%%%%%%%%%%%%%%%%%%%%%%%%%%%%

In this section, we state and prove our Main Theorem and some of its
implications.

%%%%%%%%%%%%%%%%%%%%%%%%%%%%%%%%%%%%%%%%%%%%%%%%%%%%%%%%%%%%%%%%%%%%%%%%%%
\subsection{Main Theorem}
\label{sec31}
%%%%%%%%%%%%%%%%%%%%%%%%%%%%%%%%%%%%%%%%%%%%%%%%%%%%%%%%%%%%%%%%%%%%%%%%%%

As already mentioned in Section \ref{sec1}, the main result of this paper
is an extension of optimality results for Delaunay triangulations
from finite sets to Delaunay sets, which are necessarily infinite.
Section \ref{sec2} provides all the technical concepts needed to give
a formal proof of the theorem that facilitates this result.
\begin{result}[Main Theorem]
  Let $X$ be a Delaunay set in $\Rspace^d$.
  \begin{enumerate}
    \item[(i)] In $\Rspace^2$, $F \in \Fcal$ implies
    $f (\Delaunay{X}) \leq f(T)$ for all uniformly bounded triangulations
    $T$ of $X$.  \Point{Theorem 4.2.}
    \item[(ii)] In $\Rspace^d$, $F \in \Gcal$ implies
    $f (\Delaunay{X}) \leq f(T)$ for all uniformly bounded triangulations
    $T$ of $X$.  \Point{Theorem 5.1.}
  \end{enumerate}
\end{result}
\proof
 Fix a uniformly bounded triangulation $T$ with parameter $q$ of $X$.
 We prove the inequalities by comparing subsets of $d$-simplices
 of $T$ and $D = \Delaunay{X}$.
 We begin with (ii).
 For every radius $\alpha$, we write $T(\alpha) \subseteq T$
 and $D(\alpha) \subseteq D$ for the sets of
 simplices contained in $\Bspace_\alpha$.
 Furthermore, we write $D'(\alpha) \subseteq D(\alpha)$ for the 
 set of simplices contained in the underlying space of $T(\alpha)$.
 Summing $F$ over the $d$-simplices in these sets, we have
 \begin{eqnarray}
   \Sigma_{T(\alpha)}
     &\!\!\!=\!\!\!&  [\Sigma_{T(\alpha)} - \Sigma_{D'(\alpha)}]
        + [\Sigma_{D'(\alpha)} - \Sigma_{D(\alpha)}]
        +  \Sigma_{D(\alpha)} .
   \label{eqn:Sigma}
 \end{eqnarray}
 The first difference on the right-hand side is non-negative by assumption
 of $F \in \Gcal$.
 We prove shortly that the second difference is bounded from above by
 a constant times $\alpha^{d-1}$.
 This implies that dividing by $V_d \alpha^d$ and taking the lower limit
 gives $f(T) \geq f(D)$, as required.
 To prove the bound for the second term,
 we assume $\alpha \geq 4q$,
 so that the Simplex Count Lemma implies that
 the number of $d$-simplices in $D(\alpha) - D(\alpha - 2q)$
 is bounded from above by a constant times $\alpha^{d-1}$.
 We get the same bound for $D(\alpha) - D'(\alpha)$
 because $T(\alpha)$ is uniformly bounded, with parameter $q$,
 so it covers all of $\Bspace_{\alpha-2q}$, which implies
 $D(\alpha-2q) \subseteq D'(\alpha)$.

 The proof of (i) is similar, except that we have to do more work
 to construct the sets of $d$-simplices, now triangles
 needed for the comparison.
 We assume $\alpha \geq 8q$ and let $T(\alpha)$ and $D(\alpha)$ be as before.
 We construct $T'(\alpha)$ by modifying $T(\alpha)$ through a sequence of
 directed flips applied to non-locally Delaunay edges.
 Iterating the directed flip, we can guarantee that all interior edges
 of $T' (\alpha)$ are locally Delaunay.
 A directed flip does not increase the size of the larger circumcircle
 (see e.g.\ \cite{Mus10}), so flipping does not take us outside
 the class of uniformly bounded triangulations.
 Importantly, the flips turn a large portion of $T(\alpha)$
 into Delaunay triangles, namely $D(\alpha-6q) \subseteq T'(\alpha)$.
 To see this, we note that every triangle $\Simplex \in T'(\alpha)$
 contained in $B_{\alpha-4q}$ is also in $D$.
 Indeed, its circumsphere is contained in $B_{\alpha-2q}$,
 which is contained in the underlying space of $T'(\alpha)$.
 If $\Simplex$ were not empty, we would have a vertex inside the
 circumcircle, which would imply an edge that is not locally Delaunay
 between this vertex and $\Simplex$, which is a contradiction;
 see also \cite{Del34} where this argument is used to prove
 part (ii) of the Delaunay Triangulation Theorem.
 Finally, the triangles of $T'(\alpha)$ contained in $\Bspace_{\alpha-4q}$
 cover $\Bspace_{\alpha-6q}$,
 which implies $D(\alpha-6q) \subseteq T'(\alpha)$, as required.
 For the comparison, we consider
 \begin{eqnarray}
   \Sigma_{T(\alpha)}  &=&  [\Sigma_{T(\alpha)} - \Sigma_{T'(\alpha)}]
                          + [\Sigma_{T'(\alpha)} - \Sigma_{D(\alpha-6q)}] \\
                       && + [\Sigma_{D(\alpha-6q)} - \Sigma_{D(\alpha)}]
                          +  \Sigma_{D(\alpha)} .
 \end{eqnarray}
 The first difference is non-negative,
 and the second and third differences are bounded from above by a constant
 times $\alpha^{d-1}$.
 Dividing by $V_2 \alpha^2$ and taking the lower limit,
 as $\alpha$ goes to infinity,
 we get $f(T) \geq f(D)$, as required.
\eop

%%%%%%%%%%%%%%%%%%%%%%%%%%%%%%%%%%%%%%%%%%%%%%%%%%%%%%%%%%%%%%%%%%%%%%%%%%
\subsection{Implications in the Plane}
\label{sec32}
%%%%%%%%%%%%%%%%%%%%%%%%%%%%%%%%%%%%%%%%%%%%%%%%%%%%%%%%%%%%%%%%%%%%%%%%%%

There are many functionals on triangles that are known to be in $\Fcal$.
Applying the Main Theorem thus gives many optimality results for
Delaunay triangulations of Delaunay sets.
\begin{result}[Corollary A]
  Let $\Simplex$ be a triangle in $\Rspace^2$, with edges of
  length $a, b, c$, let $c_1 > 0$ and $c_2 \geq 1$ be constants,
  and consider the following list of functionals:
  \begin{itemize}
    \item $F_1 (\Simplex) = \Circumradius{c_1}{\Simplex}$.
    \item $F_2 (\Simplex) = \Circumradius{c_2}{\Simplex} \cdot \Area{\Simplex}$.
    \item $F_3 (\Simplex) = - \Inradius{\Simplex}$.
    \item $F_4 (\Simplex) = (a^2 + b^2 + c^2) / \Area{\Simplex}$.
    \item $F_5 (\Simplex) = (a^2 + b^2 + c^2) \cdot \Area{\Simplex}$.
    \item $F_6 (\Simplex) =
                 \Edist{\Centroid{\Simplex}}{\Circumcenter{\Simplex}}^2
                 \cdot \Area{\Simplex}$.
  \end{itemize}
  Then $f_i (\Delaunay{X}) \leq f_i(T)$ for every Delaunay set
  $X \subseteq \Rspace^2$, for every uniformly bounded triangulation
  $T$ of $X$, and for $1 \leq i \leq 6$.
  \Point{Corollary 4.3.}
\end{result}
\proof
 It is easy to see that all listed functionals belong to $\Ecal$.
 For finite sets, the optimality of the Delaunay triangulation
 for $f_1$ and $f_4$ was proved in \cite{Mus97},
 for $f_2$ and $f_6$ it was proved in \cite{Mus10},
 for $f_3$ it was proved in \cite{Lam94},
 and for $f_5$ it was proved in \cite{Raj94}.
 It follows that $F_i \in \Fcal$ for $1 \leq i \leq 6$,
 so the claim follows from (i) in the Main Theorem.
\eop

%%%%%%%%%%%%%%%%%%%%%%%%%%%%%%%%%%%%%%%%%%%%%%%%%%%%%%%%%%%%%%%%%%%%%%%%%%
\subsection{Implication in $d$ Dimensions}
\label{sec33}
%%%%%%%%%%%%%%%%%%%%%%%%%%%%%%%%%%%%%%%%%%%%%%%%%%%%%%%%%%%%%%%%%%%%%%%%%%

We have one example of a functional on $d$-simplices that is in $\Gcal$,
namely the extension of $F_5$ to three and higher dimensions.
Writing $a_1$ to $a_k$ for the lengths of the
$k = {d+1 \choose 2}$ edges of a $d$-simplex $\Simplex$,
we define $F_R (\Simplex) = \Volume{\Simplex} \sum_i a_i^2$;
see also \cite{Ako09}.
Rajan proved that for finite sets in $\Rspace^d$,
the density of $F_R$ attains its minimum for the Delaunay triangulation.
We will extend his proof to show that $F_R$ belongs to $\Gcal$.
With this, we get another consequence of the Main Theorem.
\begin{result}[Corollary B]
  We have $f_R (\Delaunay{X}) \leq f_R(T)$ for every Delaunay set
  $X \subseteq \Rspace^d$ and for every uniformly bounded triangulation
  $T$ of $X$.
  \Point{Theorem 5.7.}
\end{result}
\proof
 The main tool in this proof is the lifting of a point $y \in \Rspace^d$
 to the point $y^+ = (y, \norm{y}^2) \in \Rspace^{d+1}$,
 an idea that goes back to Voronoi \cite{Vor08}.
 Note that $y^+$ lies on the graph of the function
 $\varpi: \Rspace^d \to \Rspace$
 defined by $\varpi (x) = \norm{x}^2$.
 For $Y \subseteq \Rspace^d$, write $Y^+$ for the corresponding
 set of lifted points, and let $\conv{Y^+}$ be its convex hull.
 Assuming $Y$ is finite and generic, $\conv{Y^+}$ is a convex polytope
 whose faces are simplices.
 We distinguish between \emph{lower} faces whose outward
 normals point down -- against the direction of the $(d+1)$-st
 coordinate axis -- and \emph{upper} faces whose outward normals point up.
 Importantly, if we project all lower faces vertically to $\Rspace^d$,
 then we obtain the Delaunay triangulation of $Y$.
 \Point{Lemma 5.2.}
 
 Consider now the functional $F_L$ that maps a $d$-simplex
 $\Simplex$ in $\Rspace^d$ to the $(d+1)$-dimensional volume between
 the convex hull of the $d+1$ lifted vertices and the graph of $\varpi$.
 More precisely, it is the volume of the portion of the vertical
 $(d+1)$-dimensional prism over $\Simplex$ that is bounded above by the
 convex hull of the lifted vertices and below by the graph of $\varpi$.
 It is not difficult to prove that $F_L$ is invariant under
 isometries of $\Rspace^d$, and to use this fact to show that $F_L \in \Ecal$.
 \Point{Lemmas 5.3 and 5.4.}
 The reason for our interest in $F_E$ is the relation
 \begin{eqnarray}
   F_R (\Simplex)  &=&  (d+1)(d+2) F_E (\Simplex)
 \end{eqnarray}
 proved in \cite{Raj94}.
 \Point{Lemma 5.6.}
 Since the two functionals differ only by a multiplicative constant,
 it follows that $F_R$ also belongs to $\Ecal$.
 It remains to prove that $F_R$ belongs to $\Gcal$, which we do by
 showing that $F_E$ belongs to $\Gcal$.
 Indeed, this should be clear from the lifting result:
 the lifted images of the $d$-simplices in the Delaunay triangulation
 are closer to the graph of $\varpi$ than those of other $d$-simplices.
 More specifically, if $T'$ is a simplicial complex with finite
 vertex set $Y$ in $\Rspace^d$, and all simplices of 
 $D' \subseteq \Delaunay{Y}$ are contained in the underlying space of $T'$,
 then the total $(d+1)$-dimensional volume we get for $T'$
 is larger than or equal to that we get for $D'$.
 But this implies $F_E \in \Gcal$,
 and therefore $F_R \in \Gcal$, as required.
 \Point{Lemma 5.5.}
\eop

%\clearpage
%%%%%%%%%%%%%%%%%%%%%%%%%%%%%%%%%%%%%%%%%%%%%%%%%%%%%%%%%%%%%%%%%%%%%%%%%%
\section{Discussion}
\label{sec4}
%%%%%%%%%%%%%%%%%%%%%%%%%%%%%%%%%%%%%%%%%%%%%%%%%%%%%%%%%%%%%%%%%%%%%%%%%%

The main contribution of this paper is an extension of optimality results
that hold for Delaunay triangulations of finite sets to Delaunay sets,
which are necessarily infinite.
In the plane, this extension holds for all functionals that improve
upon flipping an edge that is not locally Delaunay.
In three and higher dimensions, we need stronger properties to prove
the extension.
It would be interesting to know whether these stronger properties
are necessary.
Specifically, is it true that $F \in \Fcal$ implies that the density of
$F$ attains its minimum at the Delaunay triangulation of a Delaunay set
in $\Rspace^d$, also for $d \geq 3$?
Similarly, are there functionals in $\Gcal$ that are not in $\Fcal$,
or is $\Fcal = \Gcal$?

%%%%%%%%%%%%%%%%%%%%%%%%%%%

\end{document}